# Chapter Three

# What Mathematical Logic Says about the Foundations of Mathematics

## Claudio Bernardi

SUMMARY My purpose is to examine some concepts of mathematical logic studied by Carlo Cellucci. Today the aim of classical mathematical logic is not to guarantee the certainty of mathematics, but I will argue that logic can help us to explain mathematical activity; the point is to discuss what and in which sense logic can "explain". For example, let us consider the basic concept of an axiomatic system: an axiomatic system can be very useful to organize, present, and clarify mathematical knowledge. And, more importantly, logic is a science with its own results: so, axiomatic systems are also interesting because we know several revealing theorems about them. Similarly, I will discuss other topics such as mathematical definitions, and some relationships between mathematical logic and computer science. I will also consider these subjects from an educational point of view: can logical concepts be useful in teaching and learning elementary mathematics?

## 1. Mathematical logic vs. the foundations of mathematics

There is no doubt that research in mathematical logic can contribute to the study of the foundations of mathematics. For instance, mathematical logic provides answers (both complete and partial) to the following questions:
- Given a precisely stated conjecture, can we be sure that eventually a good enough mathematician will be able to prove or disprove it?
- Can all mathematics be formalized?

- Is there a "right" set of axioms for arithmetic, or for mathematical analysis? is there a proof (in some fixed theory) for any statement of arithmetic which is true in **N**?
- Can we prove the consistency of standard mathematical theories? and what does it mean to prove consistency?
- By adding a new axiom to a theory, we find new theorems; but can we also expect to find shorter proofs for old theorems?
- Will we ever construct a computer that will be capable of answering all mathematical problems?
- Is any function from **N** to **N** computable by an appropriate computer? if not, how can we describe computable functions?
- If we know that a computation ends, can we estimate the time necessary to complete the computation?
- Is it true that, if a "short" statement is a theorem, then there is a short proof for it?

The list could be much longer. In some cases (as in the first question) the answer given by logic contradicts the expectations of a working mathematician, while in other cases (as in the last question) the answer confirms that expectation.

However, it is not true that the general purpose of mathematical logic is to clarify the foundations of mathematics. First of all, for the past few decades, much research in logic has been of mainly technical value and does not deal directly with the foundation of mathematics. Perhaps in the nineteenth century, logic was regarded as a way to guarantee the *certainty of mathematics*. But nowadays we do not expect that much: it seems naïve, and perhaps even futile, to hope for a definitive, proven certainty of mathematics.

Let us start from the beginning. Mathematical logic provides us with a precise definition of a proof and suggests rigorous methods and procedures for developing mathematical theories. But these are just the initial steps of mathematical logic: if logic consisted only in giving detailed definitions of proofs and theories, it would not be of great scientific importance. While succeeding in formalizing statements and arguments is interesting, the historical and cultural importance of proof theory, model theory, and recursion theory strongly depends on the results achieved in these areas (for example, on the answers given to the previous questions).

In other words, mathematical logic is a way of organizing mathematics and solving paradoxes; but I find that logic is interesting also because its organization of mathematics provides significant results. In fact, any theory grows if and when results are found.

So, we can distinguish between two kinds of logical results which can be useful in the study of foundations and, more generally, to working mathematicians.

On the one hand, mathematical logic provides explicit rules that mathematicians habitually use (often without being fully aware of it), inserting them into a clear and consistent framework; in this way more complex situations can be tackled. For instance in logic:

- it is explained what a proof by contradiction, or a counterexample, is; it is not impossible for a mathematician, who in his work usually gives proofs by contradiction and counterexamples, to be unable to give clear answers to explain these totally elementary concepts[1];
- various forms of the principle of induction are stated explicitly and compared;
- equivalents to the axiom of choice, or weaker forms of it, are stated and recognized.

On the other hand, mathematical logic allows us to construct a theoretical framework that clarifies the meaning and limitations of mathematical activity. The study of logic can provide information of the following kind: this theory is decidable, while this other one is not. Note that often logical results contrast with the naïve expectations of working mathematicians.

## 2. Mathematical experience vs. mathematical logic

In the introduction to the book *The Mathematical Experience* (Davis and Hersh 1981), Gian-Carlo Rota, challenging the idea that mathematics consists mainly of the demonstration of theorems, wrote a famous sentence: "A mathematician's work is mostly a tangle of guesswork, analogy, wishful thinking and frustration, and proof, far from being the

---

[1] Let us briefly recall that the word counterexample denotes an example that shows that a statement is not correct, when the statement consists of an implication that is preceded by a universal quantifier. To this end, we have to construct an object x that satisfies the hypothesis but not the thesis. In formal terms, the explanation is clearer: to conclude that $\forall x\, [P(x) \rightarrow Q(x)]$ does not hold, we have to prove $\neg\, \forall x\, [P(x) \rightarrow Q(x)]$, that is, $\exists x\, \neg\, [P(x) \rightarrow Q(x)]$ and this formula is equivalent to $\exists x\, [P(x) \wedge \neg\, Q(x)]$.

core of discovery, is more often than not a way of making sure our minds are not playing tricks."

Of course Gian-Carlo Rota was right[2]. It is true that in mathematical experience, when checking a method, testing a tool, or hoping that an application will follow, there are very often trials and failures. But in a logic book we will not find a chapter about hope or failure: it is not the purpose of mathematical logic to describe how a mathematician works.

Regarding mathematical activity and its formalization, it is worth quoting three amusing conversations. In Davis and Hersh (1981) a philosophy grad student asks the Ideal Mathematician, "What is a mathematical proof, really?". A similar question is considered in Devlin (1992); while in Hersh (2011) a Successful Mathematician is accosted by the Stubborn Student, who has trouble when comparing the mathematical concept of a limit with its concrete applications.

But there is a different question: what is the meaning of logical notions such as mathematical proof and mathematical theories? Of course, these notions have to do with the work of a mathematician.

First of all, a distinction must be made between the way in which a mathematician works and the final presentation of a subject. Logic refers to the organization of a mathematical subject in a clear form, so that other people can understand it: the starting point (axioms), new concepts (definitions), properties and consequences (theorems).
In any case, mathematical logic simply supplies *a suitable frame* for mathematical theories, but we could also find other models. Is the logical model a convenient one? When discussing this point, we have to face two different questions:

- is the model faithful?
- is the model useful?

Let me give a rough example: in our framework, taking a photograph could provide a faithful description, but would be totally useless. A faithful description which does not yield results is much less interesting than an unfaithful description which yields results. There is no doubt that the logical formalization of the concept of a proof is far from a concrete proof and even further from the way in which a proof is found, but many results can be found in proof theory, which concern mathematicians, their expectations, and the limits of mathematics. From this point of view, *logic is a fruitful model*. We could make similar remarks about other areas of mathematics. Consider mathematical analysis: it rests on the set of real numbers, which can hardly be regarded as a faithful model of reality. But mathematical analysis has enormous importance, because its theorems can be fruitfully applied in physics and the natural sciences, in the study of the real world.

---

[2] However, I do not agree with the idea that a proof is just a method of confirming what we already know. Note that many mathematicians (including Gian-Carlo Rota) have sought elegant proofs; moreover, very often a proof allows for a deeper understanding of the subject.

Mathematical logic, exactly like mathematical analysis, is justified by its results.

## 3. Axiomatic systems

The axiomatic method is *a way of thinking*. This is true in mathematical research, as well as in mathematical education (I will deal with mathematical education in § 5 and 6). The concept of a mathematical axiom, and its meaning, is part of our culture; several sciences other than mathematics have tried to introduce something similar to axioms, in order to achieve a more convincing approach to the matter studied.

In mathematics, the axiomatic method has had enormous influence from the time of Euclid (consider the fifth postulate and non-Euclidean geometries). Or, to give a present-day example, consider *reverse mathematics* which would not have been born without the concept of an axiomatic system; see (Marcone 2009) for a general introduction.

Instead of discussing the relationship between axioms and theorems in general terms, I prefer to stress one particular aspect.

Using the terminology of recursion theory (see Rogers 1967 for definitions), the set of axioms of a typical mathematical theory is *recursive* (or *decidable*), in the sense that one can recognize if a given sentence is an axiom. On the other hand, the set of theorems is not recursive, it is only *recursively enumerable*, because, given a sentence, we are generally not able to decide whether it is a theorem or not, we can only try and hope to get an answer.

Mathematics allows us to jump from a recursive set up to a set which is recursively enumerable but not recursive. So, even if the content of any given theorem is implicit in the axioms, the set of theorems is more complex, in a technical sense, than the set of axioms.

In my opinion, and from an abstract point of view, this is the ultimate task of mathematics and, on the other hand, it explains the difficulty of doing mathematics. Some remarks are necessary.

First of all, very often in the work of a mathematician, the set of axioms is not explicitly stated; but, even if this is the case, I think it can be assumed to be recursive, in the sense that the mathematician is able to recognize the axioms he is entitled to use. Moreover, a mathematician sometimes changes his hypothesis when trying to give a proof: he hoped that a result $\alpha$ held in general in a given theory $T$, but realized that it is necessary to add a hypothesis $\beta$. The axiomatic approach is not affected in this case, because it has been proven that $\beta \rightarrow \alpha$ is a theorem of $T$.

It should be observed that the mathematical community searches continuously for new axioms, which are deeper or more general, or more suitable for some purpose, trying to give a more comprehensive explanation of a subject.

Taking account of these situations, Carlo Cellucci (2002) introduced and studied "*open systems*", which can better describe the analytical method adopted by mathematicians. A fairly similar approach was suggested in strictly logical terms by Roberto Magari (1974) and (1975). These ideas are profound and relevant, but I think that open systems are not to be juxtaposed on axiomatic systems, because both reflect mathematical activities. Note also that a new axiom, or a new hypothesis, can be considered in different situations: as a partial step in the solution of a problem, or because it allows us to find desired consequences, or when introducing a new axiomatic system.

Let us go back to the set of theorems. A working mathematician may be astonished by the fact that this set is recursively enumerable, since, at first glance, it seems to be even more complex. The point is that a mathematician does not proceed by listing the set of theorems and looking for something interesting: his procedure is completely different. However, the fact that the set of theorems is recursively enumerable corresponds to the possibility of checking any proof.

## 4. Some remarks about definitions

*a*) A definition is just an abbreviation, since a long expression is substituted by a single new term: it is introduced simply for the reader's convenience. We could always replace the defined term with the defining expression: we get a statement that is less clear than the original one, but that has the same content.

*b*) Definitions are an indispensable part of any book on mathematics. Definitions not only draw attention to what will be useful later, but correspond to concepts, and therefore must be fully understood by anybody who studies a subject. For instance, the definitions of limit, continuity, derivative, are essential for studying calculus.

In my opinion, claims *a*) and *b*) are both correct. Claim *a*) reflects a theoretical and abstract point of view, whereas *b*) applies to every human being who learns mathematics.

There is no contradiction between *a*) and *b*), just as there is no contradiction between someone who says that a sphere is a locus of points which have the same distance from a given point, and someone who thinks of a sphere as something round, that can roll perfectly in all directions.

A typical question about definitions is the following: is a theorem proved starting only from axioms or starting from both axioms and definitions?

First note that, in mathematics, a definition gives no information about the objects involved: in elementary geometry we can define the bisector of an angle before knowing that any angle has a bisector. From this point of view, *a definition cannot be a starting point for proving* something (even though, of course, it can allow us to give other definitions).

On the other hand, in many cases it seems that, during a proof, we rely on definitions, especially when basic concepts are involved. But let us compare the situation to an algebraic one. When making a complicated algebraic computation, it often happens that is convenient to introduce a new letter, for instance setting $y = x^2$ (where the letter $y$ did not occur previously). In doing this, we may be putting ourselves in a position to directly apply a known formula or to recognize an algebraic pattern; but note that the previous equivalence gives no information about $x$. The equivalence $y = x^2$ is nothing but a definition, introduced only to make what follows accessible and clear.

Returning to theorems and definitions, we conclude that any proof is supported by axioms (and already known theorems), and not by any definitions. When proving a theorem we often read statements such as: 'recalling the definition of a limit, we can say that ...'; or, '*d* is a divisor of *p*, but *p* is prime and, by the definition of a prime number, we conclude $d = 1$ or $d = p$'. However, in these cases, the reference to a definition is useful only for recalling the meaning of a term, that is, for substituting a term with the proprieties used in defining that term.

# 5. Is mathematical logic useful at school? which concepts should be taught?

Knowing mathematical logic does not provide a method for those who want to do research in geometry or analysis; in the same way, the study of logic does not provide a necessary introduction to high school mathematics. In fact, were some logical concepts to be summarized in the first chapter of a mathematical textbook, this chapter would very likely be nearly useless to students, and would soon be forgotten.

What can be useful in high school are frequent discussions of the logical aspects of the mathematical concepts and procedures the students are dealing with. Indeed, mathematical education has an educational value that does not depend on applications. I am referring to skills involving the use of language and argumentation; and language and argumentation are obviously related to logic.

Let us examine some points.

a. *Axiomatic systems*

An axiomatic system is a way both of correctly organizing and presenting knowledge in a mathematical area, and also in general of teaching a correct way of proving and deducing.

In fact, without the concept of an axiomatic system, the teaching of mathematics consists only of an unjustified set of rules and arguments, based on common sense, or on the authority of the teacher. Such teaching is appropriate to middle school, but is not always suitable for high school.

It has been remarked that, in the *Elements* of Euclid, the connections between axioms and theorems are complex: if we try to specify, for any theorem, the previous theorems and axioms upon which it depends, we find an intricate structure. Even a good student cannot completely master this logical structure.

This may be true, but I believe it is not a good reason to give up! On the contrary, the teacher will pay attention to the ties between axioms and theorems in some specific cases, stress the fact that a theorem depends or does not depend e.g. on the Pythagorean theorem, and so on, even if he knows in advance that no student will learn the entire structure of the axiomatic system of Euclidean geometry.

Equally, a comparison between Euclidean geometry and non-Euclidean geometries can be useful. The teacher will show that, in the new geometries, some known theorems remain valid, while others (such as the Pythagorean theorem, the sum of angles of a triangle, etc.) no longer hold; on the other hand, there are also new theorems (such as the fourth criterion for triangles: if the angles of a triangle are congruent to the angles of a second triangle, the two triangle are congruent to each other). We must not underestimate the educational importance of the fact that a mathematical result holds in one theory but does not hold in another.

b. *Proving*

The difference between verifying and proving is obviously fundamental when teaching and learning mathematics (incidentally, it has been noted by several people that the use of computers in geometry and arithmetic does not help in this regard: why do we need to prove what is said or shown by a computer?).

Of course, I am not referring to formal proofs. But even a student at the beginning of their high school education can understand some basic logical facts about proofs, such as:

- any proof consists of several elementary steps;
- it is not easy to find a proof, but a given proof can be checked by everybody who knows the concepts involved;
- in many cases there are different proofs for the same statement;
- some theorems are proved by contradiction;
- if a theorem is an implication, the inverse implication is not always a theorem;
- proving $\alpha \to \beta$ is logically equivalent to proving $\neg\beta \to \neg\alpha$;
- 'not for every $x$ ...' is different from 'for every $x$ not ...'.

The importance of proving in teaching and learning mathematics has been widely investigated; see for instance (Arzarello 2012), (Bernardi 1998), (Bernardi 2010), (Francini 2010).

c. *Formalizing statements*

For students at the end of high school or the beginning of university, formalizing statements is a useful exercise, and in particular finding the explicit quantifiers and implications hidden in natural language. Think of a trivial sentence like 'the square of an odd number is odd'; in formalizing it, we have to write a universal quantifier as well as an implication: for any number $n$, if $n$ is odd then $n^2$ is odd.

On the other hand, I think that translating mathematical statements into a first-order language is in general too difficult for students. Formalization can only be done in specific cases. I refer for instance to the definition of the limit of a function; in my opinion the difficulty in understanding this notion also depends on the logical structure $\forall \varepsilon\ \exists \delta\ \forall x$ (this is one of the first times a student meets three alternating quantifiers).

Many current words and expressions used in mathematics are perhaps useful in practice (and in particular in the practice of teaching), but are subject to criticism for being unclear or ambiguous: e.g., "*given* a number", "*take* a function", "we *impose* that ...", "the *general term* of a sequence", "*fix the value* $x_0$ of a variable *x*", ... Formal language allows us to clarify these situations. In this respect, logic can contribute (and in fact has contributed) to improving rigour in natural mathematical language. For other remarks, see (Bernardi 2012).

## 6. Formal notation and self-confidence at school, in algebra, geometry, logic

*Algebraic manipulation* increases self-confidence in high school students.

Of course, not all students acquire good skills in algebra; but the doubtful student seeks comfort in algebraic calculus and, in any case, tries to perform some algebraic computations. Algebraic language is effective, rules in algebra are clear and simple. Steps in computation do not require too much thought; for this very reason, students usually prefer algebra to geometry.

Students, as well as teachers, rely on algebraic formalism: this is because algebraic formalism is artificial, and therefore governed by simple rules. Sometimes the abstract is more natural than the concrete: algebraic language is abstract, in the sense that it has been built by us and for us.

From an educational point of view, there is the obvious risk of a mechanical and unconscious application of rules. On the other hand, regular practice with algebraic computation develops other skills, such as the capacity to make indirect controls: e.g., in some situations we will automatically be expecting a polynomial to be homogeneous (and therefore realize something is wrong if it turns out not to be).

The situation is completely different in the teaching and learning of *geometry*. In my opinion, the teaching difficulties of Euclidean geometry also result from the lack of convenient notation. Take for example the angles $\angle RSP$ and $\angle SPQ$ of the quadrilateral *PQRS*: while these symbols are not long, they have to be continuously interpreted within a diagram. There is no formal manipulation, with the exception of very particular cases (for instance, referring to the sum of vectors the equivalence $AB + BC = AC$ holds, but a similar equality does not hold for segments).
Unlike algebraic expressions, geometrical symbols are not suitable for computation, they have only schematic and mnemonic value.

So, in the teaching of mathematics, even though the logical structure (axioms, theorems, primitive concepts, ...) is more evident in geometry than in algebra, a formal calculus occurs only in algebra.

When introduced to the symbolic use of *connectives* and *quantifiers* at the end of high school or at the beginning of university, students are amazed by the expressive power of logical language: any mathematical statement seems to be expressible in symbolic logic.

However, at that level a *calculus* in logic can be presented very rarely, and only in some specific cases. Symbols like connectives and quantifiers allow us to express sentences in a clear and concise way, but they must be interpreted each time, and cannot be directly manipulated. This situation presents analogies with notation in geometry, rather than with algebraic symbols.

In high school we can have a "logical calculus" only in the construction of truth tables, but I think that this construction has limited value in mathematical education. From this point of view, logical language cannot provide a student with the same self-confidence as algebraic manipulation.

## 7. Mathematical logic and computer science

As is well known, computers can contribute to providing mathematical proofs. Does mathematical logic play a role in this contribution?
To answer this question, I think that a distinction has to be made: there are (at least) three different kinds of contributions made by computer science to mathematics.

- First of all, there are *computer-assisted proofs*. This name usually refers to a computer used to perform very long and complex computations, or to examine a great number of possibilities. The first famous case occurred in 1976, when the four colour theorem was proven. That proof has been discussed for a long time (perhaps too long); I believe that the point of discussion was not the certainty of the proof, but the change in the style of proving (like an athletic world record which is achieved by modifying technical equipment). A computer can also be used to find approximate solutions of equations, to simulate the evolution of a phenomenon, to visualize patterns. But these applications are not too different from the previous one. In all these cases, the role of logic is usually limited.
- The role of logic is greater when a computer is used directly to find new theorems, or new proofs of known theorems, as happens particularly in algebra and geometry. This area is related to artificial intelligence: a computer searches for new facts by combining known facts. The possibility that a computer may find and prove a theorem is fascinating.

In these cases, to plan a *theorem prover*, formalization of arguments and proofs is obviously necessary. A problem is to direct the research towards statements of some interest. It is very easy to give rules to obtain new theorems: for instance, starting from $A$ we can deduce $A \wedge A$ and $B \rightarrow A$, but how can we recognize interesting theorems among trivial ones?

- Lastly, we must mention a more recent application, where the role of formal logic is even greater.

Often, in the history of mathematics, a mistake has been found in a proof that had previously been accepted (the history of the proof of Jordan's theorem about simple closed curves in the plane is a typical example). In this third application of computer science, known as *automated proof checking* (or also, with a slightly different meaning, *automated theorem proving*), computers are used simply to formalize known proofs of theorems and check them in detail. The purpose is *to certify theorems* and collect them in libraries. We could mention, for instance, the *Mizar system* (see http://www.mizar.org/) and the *Coq Proof Assistant* (see http://coq.inria.fr).

It is currently hard to predict just how much and how widely these supports will be used, but the interesting point is that automated proof checking not only guarantees greater accuracy (even if, in any case, we cannot hope for the complete certainty of a statement): in fact, looking for a better way to formalize also sheds new light, suggests new ideas, and brings in new generalizations.